\documentclass[12pt]{amsart}
\usepackage{amscd}
\def\NZQ{\Bbb}               % the font for N,Z,Q,R,C

\def\ZZ{{\NZQ Z}}

\def\frk{\frak}               % font for "Fraktur"

\def\Phi{{\frk n}}
\def\Phi{{\frk N}}
\def\opn#1#2{\def#1{\operatorname{#2}}} % to make operators
\opn\chara{char} \opn\length{\ell} \opn\pd{pd} \opn\rk{rk}
\opn\projdim{proj\,dim} \opn\injdim{inj\,dim} \opn\rank{rank}
\opn\depth{depth} \opn\grade{grade} \opn\height{height}
\opn\embdim{emb\,dim} \opn\codim{codim}

\opn\Tr{Tr} \opn\bigrank{big\,rank}
\opn\superheight{superheight}\opn\lcm{lcm}
\opn\trdeg{tr\,deg}%\emph{
\opn\reg{reg} \opn\lreg{lreg} \opn\ini{in} \opn\lpd{lpd}
\opn\size{size}
\opn\div{div} \opn\Div{Div} \opn\cl{cl} \opn\Cl{Cl}
\opn\Spec{Spec} \opn\Supp{Supp} \opn\supp{supp} \opn\Sing{Sing}
\opn\Ass{Ass} \opn\Min{Min}
\opn\Ann{Ann} \opn\Rad{Rad} \opn\Soc{Soc}
\opn\Im{Im} \opn\Ker{Ker} \opn\Coker{Coker} \opn\Am{Am}
\opn\Hom{Hom} \opn\Tor{Tor} \opn\Ext{Ext} \opn\End{End}
\opn\Aut{Aut} \opn\id{id}

\opn\nat{nat}
\opn\pff{pf}%   \pf exists already
\opn\Pf{Pf} \opn\GL{GL} \opn\SL{SL} \opn\mod{mod} \opn\ord{ord}
\opn\Gin{Gin} \opn\Hilb{Hilb}
\opn\aff{aff} \opn\con{conv} \opn\relint{relint} \opn\st{st}
\opn\lk{lk} \opn\cn{cn} \opn\core{core} \opn\vol{vol}
\opn\link{link} \opn\star{star}
\opn\gr{gr}

\def\pot#1#2{#1[\kern-0.28ex[#2]\kern-0.28ex]}

\opn\dirlim{\underrightarrow{\lim}}
\opn\inivlim{\underleftarrow{\lim}}
\def\Implies{\ifmmode\Longrightarrow \else
        \unskip${}\Longrightarrow{}$\ignorespaces\fi}
\def\implies{\ifmmode\Rightarrow \else
        \unskip${}\Rightarrow{}$\ignorespaces\fi}
\def\iff{\ifmmode\Longleftrightarrow \else
        \unskip${}\Longleftrightarrow{}$\ignorespaces\fi}

\let\:=\colon
\newtheorem{Theorem}{Theorem}[section]
\newtheorem{Lemma}[Theorem]{Lemma}
\newtheorem{Corollary}[Theorem]{Corollary}

\newtheorem{Remark}[Theorem]{Remark}

\newtheorem{Example}[Theorem]{Example}

\let\epsilon\varepsilon
\let\phi=\varphi
\let\kappa=\varkappa
\def\qed{\ifhmode\textqed\fi
      \ifmmode\ifinner\quad\qedsymbol\else\dispqed\fi\fi}
\def\textqed{\unskip\nobreak\penalty50
       \hskip2em\hbox{}\nobreak\hfil\qedsymbol
       \parfillskip=0pt \finalhyphendemerits=0}
\def\dispqed{\rlap{\qquad\qedsymbol}}

\opn\dis{dis}
\def\pnt{{\raise0.5mm\hbox{\large\bf.}}}

\opn\Lex{Lex}

\begin{document}

\title{Sequentially Cohen-Macaulay monomial ideals of embedding dimension four}

\author{  Sarfraz Ahmad, Dorin Popescu}

\thanks{The authors are highly grateful to the School of Mathematical Sciences,
  GC University, Lahore, Pakistan in supporting and facilitating this
  research. The second author was supported by  CNCSIS and the Contract 2-CEX06-11-20/2006
 of the Romanian Ministery of Education and Research and the Higher Education Comission of Pakistan. }

\address{Sarfraz Ahmad, School of Mathematical Sciences, 68-B New Muslim Town,
    Lahore,Pakistan.}\email{sarfraz11@gmail.com}
\address{Dorin Popescu, Institute of Mathematics "Simion Stoilow",
University of Bucharest, P.O.Box 1-764, Bucharest 014700, Romania}
\email{dorin.popescu@imar.ro} \maketitle

\begin{abstract}
Let $I$ be a monomial ideal of the polynomial ring
$S=K[x_1,\ldots,x_4]$ over a field $K$. Then $S/I$   is sequentially
Cohen-Macaulay if and only if $S/I$ is pretty clean. In particular,
if $S/I$ is  sequentially Cohen-Macaulay then $I$ is a Stanley
ideal. \vskip 0.4 true cm
 \noindent
\end{abstract}

\section*{Introduction}

Let $S=K[x_1,\ldots,x_n]$ be a polynomial ring over a field $K$ and
$I\subset S$ a monomial ideal. If  $S/I$ is Gorenstein of
codimension three then a description of $I$ is given in
\cite[Theorem 6.1]{BH1} in terms of the minimal system of monomial
generators. Here we are interested to describe monomial ideals $I$
when $n=4$ and $S/I$ is Cohen-Macaulay of codimension two in terms
of the primary decomposition of $I$. As a consequence we get a
particular form of \cite[Proposition 1.4]{HJY} for $n=4$, which says
that if $S/I$ is Cohen-Macaulay  of codimension two then $S/I$ is
{\em clean}, that is (after \cite{Dr}) there exists a prime
filtration $I=F_0\subset F_1\subset\ldots \subset F_r=S$ of monomial
ideals such that $F_i/F_{i-1}\cong (S/P_i)(a_i)$ for some prime
ideals $P_i$ of $S$ with $ht(P_i)=\dim(S/I)$ and $a_i\in \ZZ$,
$i=1,\ldots,r$.

More general, given a monomial ideal $I$ of $S$ then $S/I$ is called
{\em pretty clean} after \cite{HP} if there exists a prime
filtration $I\subset F_1\subset\ldots \subset F_r=S$ of monomial
ideals such that $F_i/F_{i-1}\cong S/P_i(a_i)$ for some prime ideals
$P_i$ of $S$ with the property that $P_i\subset P_j$ and $i\leq j$
implies $P_i=P_j$, that is, roughly speaking, "bigger primes come
first" in the filtration. \cite[Corollary 4.3]{HP} says that if
$S/I$ is pretty clean then $S/I$ is sequentially Cohen-Macaulay,
that is the non-zero factors of the dimension filtration of
\cite{Sc} (see next section) are Cohen-Macaulay. Our Theorem
\ref{main} says that for $n=4$ it is true also the converse, namely
that if $S/I$ is  sequentially Cohen-Macaulay then $S/I$ is pretty
clean.

A decomposition of $S/I$ as a direct sum of linear $K$-spaces of the
form $S/I=\oplus_{i=1}^r u_iK[Z_i]$, where $u_i$ are monomials of
$S$ and $Z_i\subset \{x_1,\ldots,x_n\}$ are subsets, is called a
{\em Stanley decomposition}. Stanley \cite{St} conjectured that
there always exists such a decomposition  such that $|Z_i|\geq
\depth(S/I)$. If Stanley conjecture holds for $S/I$ then $I$ is
called a {\em Stanley ideal}. Our Corollary \ref{cor} says that if
$n=4$, $I$ is monomial and $S/I$ is sequentially Cohen-Macaulay then
$I$ is a Stanley ideal (this follows because $I$ is a Stanley ideal
whenever $S/I$ is pretty clean as says \cite[Theorem 6.5]{HP}).

\section{Sequentially Cohen-Macaulay monomial ideals of embedding dimension
four are pretty clean}

Let $S=K[x_1,\ldots,x_n]$ be a polynomial ring over a field $K$. The
following result  \cite[Proposition 1.4]{HJY} is essential in this
section.

\begin{Theorem}[Herzog-Soleyman Jahan-Yassemi] \label{hjy}
Let $I\subset S$ be a monomial ideal of height two such that $S/I$
is Cohen-Macaulay. Then $S/I$ is clean.
\end{Theorem}

The proof of Herzog, Soleyman Jahan and Yassemi  passes the problem
to the polarization, where they could use strong tools from
simplicial complex theory. In the next section we give a  direct
proof in the case $n=4$, which uses just elementary theory of
monomial ideals. With this occasion we give also a complete
description of all monomials ideals $I$ of height 2 in the case
$n=4$ with $S/I$ Cohen-Macaulay. The conditions given in this
description are sometimes difficult but they could easily give nice
examples of monomial ideals $I$ with $S/I$ not Cohen-Macaulay, but
having all associated primes of height 2 and with $S/\sqrt I$
Cohen-Macaulay (see Example \ref{ex}). Certainly if $S/I$ is
Cohen-Macaulay then $S/\sqrt I$ is too by \cite[Theorem 2.6]{HTT}
(this holds only for monomial ideals). We mention that special
descriptions of some monomial Cohen-Macaulay ideals of codimension 2
are given in \cite[Theorem 3.2]{HTT}.

Let $I\subset S$ be a monomial ideal and $I=\bigcap_{p\in
\Ass(S/I)}\ P_p$, $\sqrt{P_p}=p$,
 an irredundant primary decomposition of $I$. Set
$D_i(I)=\bigcap_{p\in \Ass^{>i}(S/I)}\ P_p$, for $-1\leq i<n$, where
$\Ass^{>i}(S/I)=\{p\in \Ass(S/I): \dim(S/p)>i\}$. We get in this way
the dimension filtration of $S/I$
$$I=D_{-1}(I)\subset D_0(I)\subset \ldots \subset D_{n-2}(I)\subset
D_{n-1}(I)=S,$$ introduced by Schenzel \cite{Sc} ($n$ is the number
of variables of $S$). $S/I$ is {\em sequentially Cohen-Macaulay} if
all non-zero factors of this filtration are Cohen-Macaulay. In the
monomial case, the notions of "sequentially Cohen-Macaulay" and
"pretty clean" are connected by the following result of
\cite[Corollary 4.3]{HP}.

\begin{Theorem}[Herzog-Popescu]\label{hp}
Let $I\subset S$ be a monomial ideal and
$$I=D_{-1}(I)\subset D_0(I)\subset \ldots \subset D_{n-2}(I)\subset
D_{n-1}(I)=S$$ the dimension filtration of $S/I$. Then the following
statements are equivalent:
\begin{enumerate}
\item $S/I$ is pretty clean,

\item $S/I$ is sequentially Cohen-Macaulay and all non-zero factors
of the dimension filtration are clean,

\item all non-zero factors
of the dimension filtration are clean.

\end{enumerate}
\end{Theorem}

From now on $S=K[x,y,z,w]$, that is the case  $n=4$. The above
theorems are main tools in proving the following:

\begin{Theorem}\label{main}
Let $I\subset S=K[x,y,z,w]$ be a monomial ideal. Then $S/I$ is
pretty clean if and only if $I$ is sequentially Cohen-Macaulay.
\end{Theorem}

\begin{proof}
By Theorem \ref{hp} it is enough to show that the
non-zero factors of the dimension filtration
$$I=D_{-1}(I)\subset D_0(I)\subset D_1(I)\subset D_2(I)\subset
D_3(I)=S,$$ are clean if they are Cohen-Macaulay. Since $S$ is
factorial ring and $D_2(I)$ is an intersection of primary height one
ideals we get $D_2(I)=(u)$ for a certain monomial $u\in S$. Clearly
$S/(u)$ is clean (see e. g. \cite[Lemma 1.9]{Ja}). As
$D_2(I)/D_1(I)\cong S/(D_1(I):u)$ is Cohen-Macaulay of dimension 2
we get $D_2(I)/D_1(I)$ clean by Theorem \ref{hjy}. Now note that
$D_1(I)/D_0(I)$ and $D_1(I)/I$ are clean by \cite[Corollary 2.2]{Po}
 because the prime ideals associated to those modules are of height
$\geq 3$.
\end{proof}

\begin{Corollary} \label{cor}
Let $I\subset S=K[x,y,z,w]$ be a monomial ideal. If $S/I$ is
sequentially Cohen-Macaulay then $I$ is a Stanley ideal.
\end{Corollary}
\begin{proof}
By the above theorem $S/I$ is pretty clean and it is enough to apply
\cite[Theorem 6.5]{HP}.
\end{proof}

\section{Proof of Theorem \ref{hjy} in the case $n=4$}

Let $K$ be a field and $S=K[x,y,z,w]$ be the polynomial ring in four
variables. We denote $G(I)$ to be the set of minimal monomial
generators for an ideal $I$ in $S$. First next lemmas , which
involve ideals generated in 3 variables are easy and contained
somehow in \cite{Ja}, but we prove them for the sake of our
completeness.

\begin{Lemma} \label{1}
Let $I\subset S$ be a monomial ideal such that
$\Ass(S/I)=\{(x,y),(x,z)\}$. Then $S/I$ is clean.
\end{Lemma}
\begin{proof} Let $I=\bigcap\limits_{i=1}^{s} Q_i$ be the
irredundant decomposition of $I$ in irreducible monomial ideals (see
\cite{Vi}). Let $Q_1=(x^a,y^b)$ and $J=\bigcap\limits_{i=2}^{s}
Q_i$, where $b$ is the maximum power of $y$, which  enters in
$G(Q_i)$. Then we have the filtration $I\subset (I,x^a)\subset Q_1
\subset S$.
\\Clearly $S/Q_1$ is clean. Apply induction on $s$. We have
$Q_1/(I,x^a)\cong S/((I,x^a):y^b)$. As $G((I,x^a):y^b)$ contains
only monomials in $\{x,z\}$ we see that $((I,x^a):y^b)$ is primary
because it is the intersection of those $(Q_i,x^a)$ with
$\sqrt{Q_i}=(x,z)$. Thus $Q_1/(I,x^a)$ is clean.
\\On the other hand $ (I,x^a)/I\cong S/(I:x^a)$ and
$(I:x^a)=\bigcap\limits_{i=2}^{s} (Q_i:x^a)$. We are done by
induction
 hypothesis on $s\geq 2$,  case
$s=2$ being trivial since $(I:x^a)$ is irreducible. Thus $(I,x^a)/I$
is clean.
\end{proof}
\begin{Lemma}
\label{2} Let $I\subset S$ be a monomial ideal such that
$$\Ass(S/I)=\{(x,y),(x,z),(y,z)\}.$$
 Then $S/I$ is clean.
\end{Lemma}
\begin{proof} Let $I=\bigcap\limits_{i=1}^{s} Q_i$ be the
irredundant decomposition of $I$ in irreducible monomial ideals. Let
$Q_1=(x^a,y^b)$ and  $J=\bigcap\limits_{i=2}^{s} Q_i$, where $b$ is
the maximum such that $y^b$ enter in $G(Q_i)$. Then we have the
filtration $I\subset (I,x^a)\subset Q_1 \subset S$.
\\Clearly $S/Q_1$ is clean. Apply induction on $s$.
 We have $Q_1/(I,x^a)\cong S/((I,x^a):y^b)$. As $G((I,x^a):y^b)$
contains only monomials in $\{x,z\}$ we see that $((I,x^a):y^b)$ is
primary and its radical is $(x,z)$. Thus $Q_1/(I,x^a)$ is clean.
\\On the other hand $S/(I:x^a)\cong (I,x^a)/I$ and $(I:x^a)=\bigcap\limits_{i=2}^{s} (Q_i:x^a)$.
 We apply  induction hypothesis on $s\geq 3$, $(I:x^a)$ being in the case $s=3$   just an irreducible ideal.
 Thus $(I,x^a)/I$ is  clean.
\end{proof} \label{3}
\begin{Lemma}
Let $I\subset S$ be a monomial ideal such that
$\Ass(S/I)=\{(x,y),(z,w)\}$. Then $S/I$ is not Cohen-Macaulay.
\end{Lemma}
\begin{proof} Let $I=P_1\cap P_2$ be the irredundant
decomposition of $I$ in monomial primary ideals, let us say
$\sqrt{P_1}=(x,y), \sqrt{P_2}=(z,w)$. Then $S/(P_1+P_2)$ has
dimension $0$ and from the exact sequence $0\rightarrow S/I
\rightarrow S/P_1\oplus S/P_2\rightarrow S/(P_1+P_2)\rightarrow 0$,
we get $depth(S/I)=1$ by Depth Lemma (see e. g. \cite[Proposition
1.2.9]{BH2}). Thus $S/I$ is not Cohen-Macaulay.
\end{proof}

\begin{Remark} {\em The above lemma is trivial when $I$ is a reduced
ideal because the simplicial complex associated to $I$ is not
connected and so not Cohen-Macaulay. If $I$ is Cohen-Macaulay then
$\sqrt I$ is too by \cite[Theorem 2.6]{HTT}, which gives another
proof of this lemma.} \end{Remark}

\begin{Lemma}\label{4}
Let $I\subset S$ be a monomial ideal such that
$$\Ass(S/I)=\{(x,y),(x,z),(x,w)\}$$ and let $I=P_1\cap P_2\cap P_3$ be
the irredundant monomial primary decomposition of $I$, where
$\sqrt{P_1}=(x,y)$, $\sqrt{P_2}=(x,z)$, $\sqrt{P_3}=(x,w)$. Then
(S/I) is clean.
\end{Lemma}
\begin{proof}
Let $I=\bigcap\limits_{i=1}^{s}Q_i$ be the irredundant monomial
irreducible decomposition of $I$. Apply induction on $s$. If $s=3$,
then $(P_i)_i$ must be irreducible and so $P_1$ has the form
$(x^a,y^b)$. We consider the filtration
\\$I\subset (I,x^a)\subset (x^a,y^b)\subset S$.
Note that $P_1/(I,x^a)\cong S/((I,x^a):y^b)$. But
$((I,x^a):y^b)=(P_1\cap (P_2,x^a)\cap ((P_3,x^a):
y^b))=((P_2,x^a):y^b)\cap ((P_3,x^a):y^b)=(P_2,x^a)\cap (P_3,x^a)$
is clean by Lemma \ref{1}. Thus $P_3/(I,x^a)$ is clean.
\\ Now note that $(I,x^a)/I\cong S/(I:x^a)$. We have $(I:x^a)=(P_2:x^a)\cap (P_3:x^a)$ and so
$S/(I:x^a)$ is clean by Lemma \ref{1}. Gluing together the clean
filtrations obtained above we get a clean filtration of $S/I$ for
$s=3$.
\\Assume $s>3$. After renumbering $Q_i$ we may suppose that
$Q_1=(x^a,y^b)$ for some $a,b$. Moreover we may suppose that $b$ is
the biggest power of $y$ which can enter in
$\bigcup\limits_{i=1}^{s} G(Q_i)$. Consider the filtration as above
$I\subset (I,x^a)\subset Q_1=(x^a,y^b)\subset S$. We have
$Q_1/(I,x^a)\cong S/((I,x^a):y^b)$ and $(I,x^a):y^b=(P_2,x^a)\cap
(P_3,x^a)$ as above. Thus $Q_1/(I,x^a)$ is clean. Now note that
$(I,x^a)/I\cong S/(I:x^a)$ and $(I:x^a)=\bigcap\limits_{i=2}^{s}
(Q_i:x^a)$ and $S/(I:x^a)$ is clean by induction hypothesis. As
above gluing the obtained clean filtrations we get $S/I$ clean.
\end{proof}
\begin{Lemma}\label{5}
Let $I\subset S$ be a monomial ideal such that
$$\Ass(S/I)=\{(x,y),(x,z),(z,w)\}$$ and let $I=P_1\cap P_2\cap P_3$ be
the irredundant monomial primary decomposition of $I$, where
$\sqrt{P_1}=(x,y)$, $\sqrt{P_2}=(x,z)$, $\sqrt{P_3}=(z,w)$. Then the
following statements are equivalent:
\\$i)$ $S/I$ is clean.
\\$ii)$ $S/I$ is Cohen-Macaulay.
\\$iii)$ $P_2\subset P_1+P_3$.
\end{Lemma}
\begin{proof}
$i)\Rightarrow ii):$\,\, By \cite[Corollary 4.3]{HP}, we get $S/I$
sequentially Cohen-Macaulay. Since all primes from $\Ass(S/I)$ have
the same dimension it follows that $S/I$ is Cohen-Macaulay.
\\$ii)\Rightarrow iii):$\,\,Let $J=P_1\cap P_2$. As in the proof of
Lemma \ref{2}, from the exact sequence $0\rightarrow S/I \rightarrow
S/J\oplus S/P_3\rightarrow S/(J+P_3)\rightarrow 0$, we get that
$depth(S/I)=1$ if $depth (S/(J+P_3))=0$. But $J+P_3=(P_1+P_3)\cap
(P_2+P_3)$ and $P_1+P_3$ is primary of height 4 and $P_2+P_3$ is
primary of height 3. Thus $depth(S/(J+P_3))=0$ if and only if
$P_2+P_3 \not\subset P_1+P_3$, that is $P_2\not\subset P_1+P_3$.
Therefore if $P_2\not\subset P_1+P_3$ then $S/I$ is not
Cohen-Macaulay, which proves $ii)\Rightarrow iii)$.
\\ $iii)\Rightarrow i):$ Suppose now that $iii)$ holds and let $I=\bigcap\limits_{i=1}^{s}
Q_i$ be the irredundant monomial irreducible decomposition of $I$.
Apply induction on $s$. If $s=3$, then $(P_i)_i$ must be irreducible
and so $P_3$ has the form $(z^r,w^t).$ We consider the filtration
\\$I\subset (I,z^r)\subset (z^r,w^t)\subset S$.
Note that $P_3/(I,z^r)\cong S/((I,z^r):w^t)$. But
$((I,z^r):w^t)=(((P_1,z^r)\cap (P_2,z^r)\cap P_3):
w^t)=((P_1,z^r):w^t)\cap ((P_2,z^r):w^t)=(P_1,z^r)\cap (P_2,z^r).$
As $P_2\subset P_1+P_3$ it follows that $P_2\subset (P_1,z^r)$ and
so $(I,z^r):w^t=(P_2,z^r)$ which is primary with
$\sqrt{(P_2,z^r)}=(x,z).$ Thus $P_3/(I,z^r)$ is clean.
 Now note
that $(I,z^r)/I\cong S/(I:z^r)$. We have $(I:z^r)=(P_1:z^r)\cap
(P_2:z^r)$ and so $S/(I:z^r)$ is clean by Lemma \ref{1}. Gluing
together the clean filtrations obtained above we get a clean
filtration of $S/I$, that is $iii) \Rightarrow i)$ for $s=3$.

Assume $s>3$. After renumbering $Q_i$ we may suppose that
$Q_1=(z^r,w^t)$ for some $r,t$. Moreover we may suppose that $t$ is
the biggest power of $w$ which can enter in
$\bigcup\limits_{i=1}^{s} G(Q_i)$. Consider the filtration as above
$I\subset (I,z^r)\subset Q_1=(z^r,w^t)\subset S$. We have
$Q_1/(I,z^r)\cong S/((I,z^r):w^t)$ and
$((I,z^r):w^t)=((P_1,z^r):w^t)\cap((P_2,z^r):w^t)\cap
((P_3,z^r):w^t)=(P_1,z^r)\cap (P_2,z^r)=(P_2,z^r)$ as above. Thus
$Q_1/(I,z^r)$ is clean. Now note that $(I,z^r)/I\cong S/(I:z^r)$ and
$(I:z^r)=\bigcap\limits_{i=2}^{s} (Q_i:z^r)$ and we apply the
induction hypothesis for $(I:z^r)$ if we see that $iii)$ holds for
it. Clearly $iii)$ implies $(P_2:z^r)\subset (P_1:z^r)+(P_3:z^r)$
which is enough (note that $(P_3:z^r)$ can be a proper ideal in this
case). As above gluing the obtained clean filtrations we get $S/I$
clean.
\end{proof}

\begin{Example} \label{ex}
{\em Let $I=(x^2,y)\cap (x,z)\cap (z,w)$. Then $S/I$ is not
Cohen-Macaulay by the above lemma, but $S/\sqrt I$ is
Cohen-Macaulay, because the simplicial complex associated to $\sqrt
I$ is shellable.}
\end{Example}

\begin{Lemma}\label{6}
Let $I\subset S$ be a monomial ideal such that
$$\Ass(S/I)=\{(x,y),(x,w),(y,w),(x,z)\}$$ and let $I=P_1\cap P_2\cap
P_3\cap P_4$ be the irredundant monomial primary decomposition of
$I$, where $\sqrt{P_1}=(x,y)$, $\sqrt{P_2}=(x,w)$,
$\sqrt{P_3}=(y,w)$, $\sqrt{P_4}=(x,z)$. Then the following
statements are equivalent:
\\$i)$ $S/I$ is clean.
\\$ii)$ $S/I$ is Cohen-Macaulay.
\\$iii)$ $P_1\subset P_3+P_4$ or $P_2\subset P_3+P_4$.
\end{Lemma}
\begin{proof}
$i)\Rightarrow ii)$ as in Lemma \ref{5}.
\\$ii)\Rightarrow iii):$\,\,Let $J=P_1\cap P_2\cap P_4$. From
 the exact sequence $0\rightarrow S/I \rightarrow S/J\oplus
S/P_3\rightarrow S/(J+P_3)\rightarrow 0$, we get that $depth(S/I)=1$
if $depth (S/J+P_3)=0$. But $(J+P_3)=(P_1+P_3)\cap (P_2+P_3)\cap
(P_4+P_3)$, where $(P_4+P_3)$ is primary of height 4 and
$(P_1+P_3),(P_2+P_3)$ are primary of height 3. Thus
$depth(S/(J+P_3))=0$ if and only if $P_1+P_3 \not\subset P_4+P_3$
and $P_2+P_3 \not\subset P_4+P_3$, that is $P_1\not\subset P_4+P_3$
and $P_2\not\subset P_4+P_3$. Therefore if $P_1\not\subset P_4+P_3$
and $P_2\not\subset P_4+P_3$ then $S/I$ is not Cohen-Macaulay, which
proves $ii)\Rightarrow iii)$.
\\ $iii)\Rightarrow i):$ Let $I=\bigcap\limits_{i=1}^{s}
Q_i$ be the irredundant monomial irreducible decomposition of $I$.
Applying induction on $s$. If $s=4$, then $(P_i)$ must be
irreducible and so $P_1$ has the form $(x^a,y^b).$ Let $iii)$ holds,
let us say $P_1\subset P_3+P_4$. Consider the filtration
\\$I\subset (I,x^a)\subset (x^a,y^b)\subset S$.
Note that $P_1/(I,x^a)\cong S/((I,x^a):y^b)$. But
$((I,x^a):y^b)=(P_1\cap (P_2,x^a)\cap (P_3,x^a)\cap ((P_4,x^a):
y^b)=((P_2,x^a):y^b)\cap ((P_3,x^a):y^b)\cap ((P_4,x^a):y^b)
=(P_2,x^a)\cap ((P_3,x^a):y^b)\cap (P_4,x^a).$\\ As $P_1\subset
P_4+P_3$ it follows that $b$ is the biggest power of $y$ appearing
in $\{G(P_1),G(P_3)\}$ and so $(I,x^a):y^b$ is generated by the
variables in $x,z,w$ only, and hence clean by Lemma \ref{1}.
\\ Now note that $(I,x^a)/I\cong S/(I:x^a)$. We have $(I:x^a)=(P_2:x^a)\cap (P_3:x^a)\cap (P_4:x^a)$,
again since by hypothesis $(I:x^a)=P_2\cap P_3$, and so $S/(I:x^a)$
is clean by Lemma \ref{2}. Gluing together the filtration described
above we get a clean filtration of $S/I$. \\ Similarly, if
$P_2\subset P_3+P_4$, and $P_2=(x^n,w^p)$, then the filtration
$I\subset (I,x^n)\subset (x^n,w^p)\subset S$ is refined  to a clean
one. That is $iii) \Rightarrow i)$ for $s=4$.
\\Assume $s>4$. After renumbering $Q_i$ we may suppose that
$Q_1=(x^a,y^b)$ for some $a,b$. Moreover we may suppose that $b$ is
the biggest power of $y$ which can enter in $ G(Q_i)$ with $\sqrt
Q_i=(x,y)$. Consider the filtration as above $I\subset
(I,x^a)\subset Q_1=(x^a,y^b)\subset S$. We have $Q_1/(I,x^a)\cong
S/((I,x^a):y^b)$ and $(I,x^a):y^b=(P_2,x^a)\cap ((P_3,x^a):y^b)\cap
(P_4,x^a)$. As $P_1\subset P_3+P_4$ it follows that $b$ is the
biggest power of $y$, which appear in $G(P_3)$. Thus
$(I,x^a):y^b=(P_2,x^a)\cap (P_4,x^a)$
 and so $Q_1/(I,x^a)$ is clean by Lemma \ref{1}. Now note
that $(I,x^a)/I\cong S/(I:x^a)$ and
$(I:x^a)=\bigcap\limits_{i=2}^{s} (Q_i:x^a)$ and we apply the
induction hypothesis for $(I:x^a)$ if we see that $iii)$ holds for
it. Clearly $iii)$ implies $(P_1:x^a)\subset (P_3:x^a)+(P_4:x^a)$
which is enough. As above gluing the described clean filtration we
get $S/I$ clean. Similarly for $P_2\subset P_3+P_4$, choosing
$Q_1=(x^n,w^p)$, we complete the proof as above.
\end{proof}
\begin{Lemma}\label{7}
Let $I\subset S$ be a monomial ideal such that
$$\Ass(S/I)=\{(x,y),(x,z),(z,w),(y,w)\}$$ and let $I=P_1\cap P_2\cap
P_3\cap P_4$ be the irredundant monomial primary decomposition of
$I$, where $\sqrt{P_1}=(x,y)$, $\sqrt{P_2}=(x,z)$,
$\sqrt{P_3}=(z,w)$, $\sqrt{P_4}=(y,w)$. Then the following
statements are equivalent:
\\$i)$ $S/I$ is clean.
\\$ii)$ $S/I$ is Cohen-Macaulay.
\\$iii)$ \{$P_1\subset P_2+P_4$ or $P_3\subset P_2+P_4$\}
          and  \{$P_2\subset P_1+P_3$ or $P_4\subset P_1+P_3$\}.
\end{Lemma}
\begin{proof}
$i)\Rightarrow ii)$ as in Lemma \ref{5}.
\\$ii)\Rightarrow iii):$\,\,Let $J=P_1\cap P_2\cap P_3$. From
the exact sequence $0\rightarrow S/I \rightarrow S/J\oplus
S/P_4\rightarrow S/(J+P_4)\rightarrow 0$, we get that $depth(S/I)=1$
if $depth (S/J+P_4)=0$. But $(J+P_4)=(P_1+P_4)\cap (P_2+P_4)\cap
(P_4+P_3)$, where $(P_2+P_4)$ is primary of height 4 and
$(P_1+P_4),(P_3+P_4)$ are primary of height 3. Thus
$depth(S/J+P_4)=0$ if and only if $P_1+P_4 \not\subset P_2+P_4$ and
$P_3+P_4 \not\subset P_2+P_4$, that is $P_1\not\subset P_2+P_4$ and
$P_3\not\subset P_2+P_4$. Therefore if $P_1\not\subset P_2+P_4$ and
$P_3\not\subset P_2+P_4$ then $S/I$ is not Cohen-Macaulay. \\ On the
other hand if $J=P_1\cap P_2\cap P_4$ then the exact sequence
$0\rightarrow S/I \rightarrow S/J\oplus S/P_3\rightarrow
S/(J+P_3)\rightarrow 0$ gives the other conditions i.e. $P_2\subset
P_1+P_3$ or $P_4\subset P_1+P_3$. Remaining choices for $J$, are
equivalent to these two cases, which proves $ii)\Rightarrow iii)$.
\\ $iii)\Rightarrow i):$ Suppose now that $iii)$ holds,  let us say
 $\{P_1\subset P_2+P_4, P_2\subset P_1+P_3\}$ holds.
 Let $I=\bigcap\limits_{i=1}^{s}
Q_i$ be the irredundant monomial irreducible decomposition of $I$.
Apply induction on $s$. If $s=4$, then $(P_i)$ must be irreducible
and so $P_1$ has the form $(x^a,y^b).$ We consider the filtration
\\$I\subset (I,x^a)\subset (x^a,y^b)\subset S$.
Note that $P_1/(I,x^a)\cong S/((I,x^a):y^b)$. But
$(I,x^a):y^b=(P_1\cap (P_2,x^a)\cap (P_3,x^a)\cap (P_4,x^a)): y^b=
((P_2,x^a):y^b)\cap ((P_3,x^a):y^b)\cap ((P_4,x^a):y^b)
=(P_2,x^a)\cap (P_3,x^a)\cap ((P_4,x^a):y^b).$\\ As $P_1\subset
P_2+P_4$, so $b$ is the biggest power of $y$ in $\{G(P_1),G(P_4)\}$.
It follows that $(I,x^a):y^b=(P_2,x^a)\cap (P_3,x^a)$. Since
$P_2\subset P_1+P_3$ it follows that $(P_2,x^a)=P_2\subset
(P_3,x^a).$ Thus $(I,x^a):y^b$ is primary and so clean.
\\ Now note that $(I,x^a)/I\cong S/(I:x^a)$. We have
$(I:x^a)=(P_2:x^a)\cap (P_3:x^a)\cap (P_4:x^a)$. As above $a$ is the
biggest power of $x$ in $G(P_2)$ because $P_1\subset P_2+P_4$. Thus
$I:x^a=P_3\cap P_4$ and so $S/(I:x^a)$ is clean by again Lemma
\ref{1}.
 Gluing together the clean filtrations obtained above we get a clean
filtration of $S/I$, that is when $s=4$, then $i)$ holds for
$\{P_1\subset P_2+P_4, P_2\subset P_1+P_3\}$.
\\Assume $s>4$. After renumbering $Q_i$ we may suppose that
$Q_1=(x^a,y^b)$ for some $a,b$. Moreover we may suppose that $b$ is
the biggest power of $y$ which can enter in $ G(Q_i)$ with $\sqrt
Q_i=(x,y)$. Consider the filtration as above $I\subset
(I,x^a)\subset Q_1=(x^a,y^b)\subset S$. We have $Q_1/(I,x^a)\cong
S/((I,x^a):y^b)$ and $(I,x^a):y^b=(P_2,x^a)\cap (P_3,x^a)\cap
((P_4,x^a):y^b)=(P_2,x^a)\cap (P_3,x^a)$ as above because
$P_1\subset P_2+P_4$. Since $P_2\subset P_1+P_3$ we see that
$(P_2,x^a)\subset (P_3,x^a)$ and so $(I,x^a):y^b=(P_2,x^a)$ is
primary. Thus $Q_1/(I,x^a)$ is clean. Now note that $(I,x^a)/I\cong
S/(I:x^a)$ and $(I:x^a)=\bigcap\limits_{i=2}^{s} (Q_i:x^a)$ and we
apply the induction hypothesis for $(I:x^a)$ if we see that
 $(P_1:x^a)\subset (P_2:x^a)+(P_4:x^a)$ and   $(P_2:x^a)\subset (P_1:x^a)+(P_3:x^a)$
 which is clear. As above gluing the obtained clean filtrations we
get $S/I$ clean.
\\Other cases from $iii)$, i.e. $\{P_1\subset P_2+P_4, P_4\subset
P_1+P_3\}$, $\{P_3\subset P_2+P_4, P_2\subset P_1+P_3\}$ and
$\{P_3\subset P_2+P_4, P_4\subset P_1+P_3\}$ are similar.
\end{proof}

\begin{Lemma} \label{8}
Let $I\subset S$ be a monomial ideal such that
$$\Ass(S/I)=\{(x,y),(x,z),(z,w),(y,w),(y,z)\}$$ and let $I=P_1\cap P_2\cap
P_3\cap P_4\cap P_5$ be the irredundant monomial primary
decomposition of $I$, where $\sqrt{P_1}=(x,y)$, $\sqrt{P_2}=(x,z)$,
$\sqrt{P_3}=(z,w)$, $\sqrt{P_4}=(y,w)$, $\sqrt{P_5}=(y,z)$. Then the
following statements are equivalent:
\\$i)$ $S/I$ is clean.
\\$ii)$ $S/I$ is Cohen-Macaulay.
\\$iii)$ \{$P_1\subset P_2+P_4$ or $P_3\subset P_2+P_4$ or $P_5\subset P_2+P_4$\}
          and  \{$P_2\subset P_1+P_3$ or $P_4\subset P_1+P_3$ or $P_5\subset P_1+P_3$\}.
\end{Lemma}
\begin{proof}
$i)\Rightarrow ii)$ as in Lemma \ref{5}.
\\$ii)\Rightarrow iii):$\,\,Let $J=P_1\cap P_2\cap P_3\cap P_5$. From
the exact sequence $0\rightarrow S/I \rightarrow S/J\oplus
S/P_4\rightarrow S/(J+P_4)\rightarrow 0$, we get that $depth(S/I)=1$
if $depth (S/J+P_4)=0$. But $(J+P_4)=(P_1+P_4)\cap (P_2+P_4)\cap
(P_3+P_4)\cap (P_5+P_4)$, where $(P_2+P_4)$ is primary of height 4
and $(P_1+P_4),(P_3+P_4),(P_4+P_5)$ are primary of height 3. Thus
$depth(S/J+P_4)=0$ if and only if $P_1+P_4 \not\subset P_2+P_4$ and
$P_3+P_4 \not\subset P_2+P_4$ and $P_3+P_4 \not\subset P_5+P_4$,
that is $P_1\not\subset P_2+P_4$ and $P_3\not\subset P_2+P_4$ and
$P_5\not\subset P_2+P_4$. Therefore if $P_1\not\subset P_2+P_4$ and
$P_3\not\subset P_2+P_4$ and $P_5\not\subset P_2+P_4$ then $S/I$ is
not Cohen-Macaulay. \\ On the other hand if $J=P_1\cap P_2\cap
P_4\cap P_5$ then the exact sequence $0\rightarrow S/I \rightarrow
S/J\oplus S/P_3\rightarrow S/(J+P_3)\rightarrow 0$ gives the other
conditions i.e. $P_2\subset P_1+P_3$ or $P_4\subset P_1+P_3$ or
$P_5\subset P_1+P_3$. Remaining choices for $J$, are equivalent to
these two cases, which proves $ii)\Rightarrow iii)$.
\\ $iii)\Rightarrow i):$ Suppose now that $iii)$ holds, let us say
 $\{P_1\subset P_2+P_4, P_2\subset P_1+P_3\}$ holds.
 Let $I=\bigcap\limits_{i=1}^{s}
Q_i$ be the irredundant monomial irreducible decomposition of $I$.
Apply induction on $s$. If $s=5$, then $(P_i)$ must be irreducible
and so $P_1$ has the form $(x^a,y^b).$
\\Here we can suppose $b$ to be the biggest power of $y$ in
$\{P_1,P_4\}$ because $P_1\subset P_2+P_4$. If $y^b\in G(P_5)$ then
we consider the filtration $I\subset (I,x^a)\subset (x^a,y^b)\subset
S$. Note that $P_1/(I,x^a)\cong S/((I,x^a):y^b)$. But
$((I,x^a):y^b)=(P_2,x^a)\cap(P_3,x^a)=(P_2,x^a)$ because $P_2\subset
P_1+P_3$. Thus $P_1/(I,x^a)$ is clean. Also note that
$(I,x^a)/I\cong S/(I:x^a)$ and $I:x^a=P_3\cap P_4\cap P_5$ because
$P_1\subset P_2+P_4$. Thus $(I,x^a)/I$ is clean by Lemma \ref{2}. If
$y^b\not \in G(P_5)$ then let $P_5=(y^r,z^t)$ and we consider the
filtration $I\subset (I,z^t)\subset (y^r,z^t)\subset S$. As above we
have $P_5/(I,z^t)\cong S/((I,z^t):y^r)$ and
$((I,z^t):y^r)=(P_2,z^t)\cap(P_3,z^t)$.  Thus $P_5/(I,z^t)$ is clean
by Lemma \ref{1}. Also note that $(I,z^t)/I\cong S/(I:z^t)$. Since
$I:z^t=P_1\cap (P_2:z^t)\cap (P_3:z^t)\cap P_4$ we see that
$(I,z^t)/I$ is clean by Lemma \ref{7}.
 Gluing together the clean filtrations obtained above we get a clean
filtration of $S/I$, that is when $s=5$, then $i)$ holds for
$\{P_1\subset P_2+P_4, P_2\subset P_1+P_3\}$.
\\Assume $s>5$. After renumbering $Q_i$ we may suppose that
$Q_1=(x^a,y^b)$ for some $a,b$. Moreover we may suppose that $b$ is
the biggest power of $y$ which can enter in $ G(Q_i)$ with $\sqrt
Q_i=(x,y)$. If $y^b\in G(P_5)$ consider the filtration as above
$I\subset (I,x^a)\subset Q_1=(x^a,y^b)\subset S$. We have
$Q_1/(I,x^a)\cong S/((I,x^a):y^b)$ and $(I,x^a):y^b=(P_2,x^a)\cap
(P_3,x^a)$ because $P_1\subset P_2+P_4$. Also we get $x^a\in
G(P_2)$. Since $P_2\subset P_1+P_3$ we have $P_2\subset (P_3,x^a)$
and so $(I,x^a):y^b=P_2$ is primary. Thus $Q_1/(I,x^a)$ is clean.
Now note that $(I,x^a)/I\cong S/(I:x^a)$ and
$(I:x^a)=\bigcap\limits_{i=2}^{s} (Q_i:x^a)$ and we apply the
induction hypothesis  because $(I:x^a)$ satisfies the condition
similar to $iii)$. Gluing the obtained clean filtrations we get
$S/I$ clean. If $y^b\not\in G(P_5)$ then $y^r\in G(P_5)$ for some
$r>b$. After renumbering $Q_i$ we may suppose that $Q_1=(y^r,z^t)$.
We consider the filtration $I\subset (I,z^t)\subset Q_1\subset S$.
We have $Q_1/(I,z^t)\cong S/((I,z^t):y^r)$ and
$((I,z^t):y^r)=(P_2,z^t)\cap (P_3,z^t)$ and applying Lemma \ref{1}
we get  $Q_1/(I,z^t)$ clean. Now the proof goes as above.
\\Other cases from $iii)$ are similar.
\end{proof}
\begin{Lemma}\label{9}
Let $I\subset S$ be a monomial ideal such that
$$\Ass(S/I)=\{(x,y),(x,z),(z,w),(y,w),(y,z),(x,w)\}$$ and let $I=P_1\cap P_2\cap
P_3\cap P_4\cap P_5\cap P_6$ be the irredundant monomial primary
decomposition of $I$, where $\sqrt{P_1}=(x,y)$, $\sqrt{P_2}=(x,z)$,
$\sqrt{P_3}=(z,w)$, $\sqrt{P_4}=(y,w)$, $\sqrt{P_5}=(y,z)$,
$\sqrt{P_6}=(x,w)$. Then the following statements are equivalent:
\\$i)$ $S/I$ is clean.
\\$ii)$ $S/I$ is Cohen-Macaulay.
\\$iii)$ \{$P_1\subset P_5+P_6$ or $P_2\subset P_5+P_6$ or $P_3\subset P_5+P_6$ or $P_4\subset P_5+P_6$\}
           \\and \{$P_1\subset P_2+P_4$ or $P_3\subset P_2+P_4$ or $P_5\subset P_2+P_4$ or $P_6\subset P_2+P_4$\}
          \\and  \{$P_2\subset P_1+P_3$ or $P_4\subset P_1+P_3$ or $P_5\subset P_1+P_3$ or $P_6\subset P_1+P_3$\}.
\end{Lemma}
\begin{proof}
$i)\Rightarrow ii)$ as in Lemma \ref{5}.
\\$ii)\Rightarrow iii):$\,\,Let $J=P_1\cap P_2\cap P_3\cap P_4\cap P_5$. From
 the exact sequence $0\rightarrow S/I \rightarrow
S/J\oplus S/P_6\rightarrow S/(J+P_6)\rightarrow 0$, we get that
$depth(S/I)=1$ if $depth (S/J+P_6)=0$. But $(J+P_6)=(P_1+P_6)\cap
(P_2+P_6)\cap (P_3+P_6)\cap (P_4+P_6)\cap (P_5+P_6)$, where
$(P_5+P_6)$ is primary of height 4 and
$\{(P_1+P_6),(P_2+P_6),(P_3+P_6),(P_4+P_6)\}$ are primary of height
3. Thus $depth(S/J+P_6)=0$ if and only if $P_1+P_6 \not\subset
P_5+P_6$ and $P_2+P_6 \not\subset P_5+P_6$ and $P_3+P_6 \not\subset
P_5+P_6$ and $P_4+P_6 \not\subset P_5+P_6$, that is $P_1\not\subset
P_5+P_6$ and $P_2\not\subset P_5+P_6$ and $P_3\not\subset P_5+P_6$
and $P_4\not\subset P_5+P_6$. So this gives one condition of $iii)$.\\
On the other hand if $J=P_1\cap P_2\cap P_3\cap P_5\cap P_6$ then
the exact sequence $0\rightarrow S/I \rightarrow S/J\oplus
S/P_4\rightarrow S/(J+P_4)\rightarrow 0$ gives the second condition
of $iii)$. And finally if $J=P_1\cap P_2\cap P_4\cap P_5\cap P_6$
then the exact sequence $0\rightarrow S/I \rightarrow S/J\oplus
S/P_3\rightarrow S/(J+P_3)\rightarrow 0$ gives the second condition
of $iii)$. Remaining choices for $J$, are equivalent to these three
cases, which proves $ii)\Rightarrow iii)$.
\\ $iii)\Rightarrow i):$ Suppose now that $iii)$ holds,  let us say
 $\{P_1\subset P_5+P_6, P_1\subset P_2+P_4, P_2\subset P_1+P_3\}$ holds.
 Let $I=\bigcap\limits_{i=1}^{s}
Q_i$ be the irredundant monomial irreducible decomposition of $I$.
Apply induction on $s$. If $s=6$, then $(P_i)$ must be irreducible
and so $P_1$ has the form $(x^a,y^b).$ We consider the filtration
\\$I\subset (I,x^a)\subset (x^a,y^b)\subset S$.
 Note that $P_1/(I,x^a)\cong
S/((I,x^a):y^b)$. But $(I,x^a):y^b
=(P_2,x^a)\cap (P_3,x^a)\cap ((P_4,x^a):y^b)\cap ((P_5,x^a):y^b)\cap (P_6,x^a).$\\
As $P_1\subset P_2+P_4$ and $P_1\subset P_5+P_6$,  $b$ is biggest
power of $y$ in $\{G(P_1),G(P_4),G(P_5)\}$ and thus
$(I,x^a):y^b=(P_2,x^a)\cap (P_3,x^a)\cap (P_6,x^a)$. Also since
$P_2\subset P_1+P_3$ it follows that $(P_2,x^a)=P_2\subset
(P_3,x^a)$. Thus $(I,x^a):y^b=(P_2,x^a)\cap (P_6,x^a)$ and
$P_1/(I,x^a)$ is clean by Lemma \ref{1}.
\\ Now note that $(I,x^a)/I\cong S/(I:x^a)$. We have
$(I:x^a)=(P_2:x^a)\cap (P_3:x^a)\cap (P_4:x^a)\cap (P_5:x^a)\cap
(P_6:x^a)$. As above $a$ is the biggest power of $x$ in
$G(P_1),G(P_2),G(P_6)$. It follows
 $I:x^a=P_3\cap P_4\cap P_5$, so $S/(I:x^a)$ is
clean by  Lemma \ref{2}.
 Gluing together the clean filtrations obtained above we get a clean
filtration of $S/I$, that is when $s=5$, then $i)$ holds for
$\{P_1\subset P_5+P_6, P_1\subset P_2+P_4, P_2\subset P_1+P_3\}$.
\\Assume $s>5$. After renumbering $Q_i$ we may suppose that
$Q_1=(x^a,y^b)$ for some $a,b$. Moreover we may suppose that $b$ is
the biggest power of $y$ which can enter in $ G(Q_i)$ with $\sqrt
Q_i=(x,y)$. Consider the filtration as above $I\subset
(I,x^a)\subset Q_1=(x^a,y^b)\subset S$. We have $Q_1/(I,x^a)\cong
S/((I,x^a):y^b)$ and $(I,x^a):y^b=(P_2,x^a)\cap
(P_3,x^a)\cap(P_6,x^a)$ because $P_1\subset P_2+P_4$, $P_1\subset
P_5+P_6$. We get also $x^a\in P_2$.
 Since $P_2\subset P_1+P_3$ we have $(P_2,x^a)\subset (P_3,x^a)$
and so $(I,x^a):y^b=(P_2,x^a)\cap (P_6,x^a)$. Thus $Q_1/(I,x^a)$  is
clean by Lemma \ref{1}.  Now note that $(I,x^a)/I\cong S/(I:x^a)$
and $(I:x^a)=\bigcap\limits_{i=2}^{s} (Q_i:x^a)$ and we apply the
induction hypothesis for $(I:x^a)$ because the condition $iii)$ are
fulfilled in this case. As above gluing the obtained clean
filtrations we get $S/I$ clean.
\\Other cases from $iii)$, are similar.
\end{proof}

\end{document}